\documentclass{article}[11pt]
\usepackage{rawfonts,times,amsfonts,graphicx}
\newtheorem{theorem}{Theorem}

\newtheorem{corollary}{Corollary}

\title{On Cycles in Random Graphs}
\author{Madhav P. Desai \\ Department of Electrical Engineering \\ Indian Institute of Technology \\ Mumbai, India} 

\begin{document}
\maketitle

\begin{abstract}
We consider the geometric 
random (GR) graph on the $d-$dimensional torus with
the $L_\sigma$ distance measure ($1 \leq \sigma \leq \infty$).
Our main result is an 
exact characterization of the probability that
a particular labeled cycle exists in this random graph. 
For $\sigma = 2$ and $\sigma = \infty$, we use 
this characterization to derive a series which evaluates to
the cycle probability.  We thus obtain 
an exact formula for the expected number of 
Hamilton cycles in the random graph (when $\sigma = \infty$
and $\sigma = 2$).  We also consider the
adjacency matrix of the random graph and derive 
a recurrence relation for the
expected values of the elementary
symmetric functions evaluated on the eigenvalues
(and thus the determinant) of the adjacency
matrix, and a recurrence relation for the expected value of 
the permanent of the adjacency matrix.   The cycle
probability features prominently in these recurrence
relations.  

We calculate these quantities for geometric random graphs
(in the $\sigma = 2$ and $\sigma = \infty$ case) with up to 
$20$ vertices, and compare them with the corresponding quantities for the 
Erd\"{o}s-R\'{e}nyi (ER) random graph with the same edge
probabilities.   The calculations indicate that the threshold for rapid
growth in the number of Hamilton cycles (as well as that for rapid
growth in the permanent of the adjacency matrix) in the GR graph
is lower than in the ER graph.  However, as the number of vertices $n$
increases, the difference between the GR and ER thresholds reduces, and
in both cases, the threshold $\sim \log(n)/n$.
Also, we observe that the expected determinant can take very large values.  
This throws some
light on the question of the maximal determinant of symmetric
$0/1$ matrices.
\end{abstract}

\section{Overview}

Consider the $d-$dimensional unit torus ${\bf T_d} \ = \ [0,1]^d$.
For $0 < r \leq 1/2$, $ 1 \leq \sigma < \infty$, the geometric random (GR) graph 
$Q_n^{(\sigma,d)}(r)$ is defined as follows.  The vertex
set corresponds to $n$ points $X_n  = \ \{ x_1, x_2, \ldots x_n \}$ 
distributed uniformly and independently 
in  $T_d$.  The set of edges $E(Q_n^{(\sigma,d)}(r))$ is defined as 
\begin{displaymath}
E(Q_n^{(\sigma,d)}(r)) \ = \ \{ \{x_i, x_j\} : \ \parallel x_i - x_j \parallel_q\   \leq r\}
\end{displaymath}
where $\parallel . \parallel_q$ is the $L_q$ norm.
Then, $Q_n^{(\sigma,d)}(r)$ is a random graph.  In this
random graph model, the presence of an edge is not 
necessarily independent of the presence of other
edges.

Another random graph model which has been very well studied is
the Erd\"{o}s-R\'{e}nyi (ER) random graph, which is defined
as follows.
Given a number $p$, $0 < p \leq 1$, let $H(n,p)$ denote the 
graph which has the vertex set $\{1,2, \ldots n\}$ and 
an edge set consisting of edges selected with probability $p$ 
(a particular edge $\{i,j\}$ is present with probability $p$ and the presence
of each edge is independent of the presence of other edges).
The ER random graph has been extensively studied. Specifically, the asymptotic
behaviour (or evolution) of this random graph has received considerable 
attention \cite{ref:ErdosRenyi, ref:Bollobas}.  The most celebrated
result of this type \cite{ref:ErdosRenyi} can be summarized as follows:
if $p = p(n) = (\log n  + c_n) /n$, then the random graph $G_n$
is almost surely connected (as $n \rightarrow \infty$) if $c_n \rightarrow \infty$,
and is almost surely disconnected if $c_n \rightarrow -\infty$).
Similar {\em thresholds} exist for all monotone graph properties\footnote{A property P is said
to be monotone if, given that it holds on a graph $G$, it also holds on $G + e$, where $e$
is an edge connecting two vertices in G.} \cite{ref:Friedgut}. 

The geometric random graph appears to exhibit similar asymptotic
properties.  In \cite{ref:GuptaKumar}, a sharp threshold for connectivity
has been exhibited for the geometric random graph on 
the unit square ($d=2$ and $\sigma=2$):
if $r = r(n)$ and if $\pi r(n)^2 = (\log n + c_n)/n$ then the random 
geometric graph is almost surely connected if $c_n \rightarrow \infty$, and
is almost surely disconnected if $c_n \rightarrow -\infty$.  The existence
of sharp thresholds for monotone properties in geometric random graphs 
has been demonstrated in \cite{ref:Goel}.  
The monograph \cite{ref:Penrose}
summarizes threshold characterizations of several connectivity related
properties of the geometric random graph.
Upper and lower bounds on the diameter of a geometric random graph in the unit
ball have been derived in \cite{ref:Ellis}. The mixing times of random walks in
geometric random graphs have been characterized in \cite{ref:Avin}. 
The limiting distribution of the eigenvalues of the adjacency matrix
of a random graph has been studied in \cite{ref:Blackwell}, \cite{ref:Rai}.   An asymptotic
bound for the second largest eigenvalue of the adjacency matrix of
a geometric random graph has been derived in \cite{ref:Boyd}.  Thus,
there is a large body of work on the asymptotic properties of
a geometric random graph.

In the finite case, one is interested in the exact formula for the
appearance of a certain property in a geometric random graph.
An example of such a characterization is an exact formula for the probability
of connectivity of a geometric random graph on a $1$-dimensional {\em unit cube} \cite{ref:DesaiManjunath},
and an exact formula for the probability of existence of a particular labeled subgraph
in the geometric random graph constructed in the $d-$dimensional unit cube using the
$L_{\infty}$ measure \cite{ref:DesaiManjunath2}.  
We will consider the finite case, and prove an exact characterization
of the probability that a labeled cycle appears in the random 
graph $Q_n^{(\sigma,d)} (r)$ (valid for $1 \leq \sigma \leq \infty$, and for all
$d \geq 1$).
Using this characterization, we show that it is possible to get exact formulas
and recurrences for the computation of quantities which are related to
cycle probabilities.  In particular, we obtain
\begin{enumerate}
\item an exact formula for the appearance of a particular labeled
cycle in $Q_n^{(\sigma,d)}(r)$ for $\sigma = 2$ and for $\sigma = \infty$ (the calculation
of the corresponding cycle probability for $H(n,p)$ is trivial, because the
edges in $H(n,p)$ are independent of each other). 
This formula
immediately yields an expression for the expected number of Hamilton cycles
in the random graph.
\item a recurrence relation for the expected values of the elementary symmetric
functions evaluated at the eigenvalues of the adjacency
matrix (as a special case, the expected value of the determinant
of the adjacency matrix) of $H(n,p)$ and $Q_n^{(\sigma,d)}(r)$.
\item a recurrence relation for the expected values of the permanent of the
adjacency matrix of $H(n,p)$ and $Q_n^{(\sigma,d)}(r)$.
\end{enumerate}

These formulas can be evaluated explicitly and provide
concrete information about random graphs with a finite
number of vertices.
For example, we observe that
cycles appear earlier in GR graphs than 
in the ER graph.  Specifically, 
the edge-probability threshold at which the expected number of Hamilton cycles  
crosses $1$ is lower in the GR graph than in the
ER graph.  However, the difference between the two thresholds
reduces as $n$ increases.  A similar observation can be
made about the expected value of the permanent. 
The expected value of
the determinant can be very different in
the GR and ER models, indicating that for particular
values of edge probabilities, the distribution of 
graphs in the GR and ER models can be very different.  
Another interesting observation is that as the edge probability is
varied between $0$ and $1$, the expected values of the determinants of
the adjacency matrix can be quite large.  In effect, these
expected values provide us some useful information about the largest
possible determinant of a symmetric $0/1$ matrix.

\section{Preliminaries}

We introduce some notation and summarize some well known
results to be used in the subsequent sections.

We use $G_n$ to denote a random graph  on $n$ vertices (in one of the models
described above).  Then $A_{G_n} = [a_{ij}({G_n})] $ is the
adjacency matrix of ${G_n}$, which is a symmetric random matrix 
with $0/1$ entries (the entries of
this matrix are correlated if $G_n$ is the GR random graph).  

Let ${\bf R}$ and ${\bf C}$ represent the
sets of real and complex numbers respectively, and let ${\bf R}^d$, ${\bf C}^d$ denote
the $d-$dimensional spaces of real and complex d-tuples.  
The  set of integers is represented by ${\bf Z}$, and
${\bf Z}^d$ is the subset of ${\bf R}^d$ consisting of 
$d-$tuples of integers.  
Elements of these spaces will be denoted by bold 
letters such as ${\bf x,y, \omega}$.  Each ${\bf x}$ in any of
these spaces is a $d-$tuple $(x_1,x_2, \ldots x_d)$.  
We will use ${\bf 1} \in {\bf Z}^d$ to denote 
the $d-$tuple with each of its entries being $1$.  
If ${\bf x} = (x_1, x_2, \ldots x_d)$ and 
${\bf y} = (y_1,y_2 \ldots y_d)$ are two elements of these
spaces, then the {\em inner product} ${\bf x . y}$ is 
$\sum x_j y_j$.
The 
$L_\sigma$ norm for these spaces defined in the usual way, 
and for ${\bf x}$, $\parallel {\bf x} \parallel_\sigma$ denotes the $L_\sigma$ norm of ${\bf x}$.
If $S \subset {\bf R}^d$, then $\Xi_S$ is the indicator 
function of $S$, so that
\begin{displaymath}
\Xi_S({\bf x}) \ = \ \left\{ \begin{array}{ll} 1 & {\rm if \ } {\bf x} \in S \\
0 & {\rm otherwise} \end{array} \right.
\end{displaymath}

For an absolutely integrable function $f:{\bf R}^d \rightarrow {\bf R}$, the
Fourier transform $\hat{f}:{\bf R}^d \rightarrow {\bf C}$ is defined as
\begin{displaymath}
\hat{f}({\bf \omega}) \ = \ \int_{{\bf x} \in {\bf R}^d} \ e^{- i {\bf \omega . x}} \ f({\bf x}) d\mu ({\bf x})
\end{displaymath}
where $d\mu({\bf x})$ is the volume element in ${\bf R}^d$ at ${\bf x}$.  Further, if $f({\bf x}) = f(-{\bf x})$
for all ${\bf x} \in {\bf R}^d$, then $\hat{f}({\bf \omega}) = \hat{f}(-{\bf \omega})$ for all
${\bf \omega} \in {\bf R}^d$, and $\hat{f}$ always takes on real values.
If $f$ is an absolutely integrable function with bounded support, and we define
\begin{equation}
f_p({\bf x}) \ = \ \sum_{{\bf u} \in {\bf Z}^d} \ f({\bf x} - {\bf u})
\end{equation}
then $f_p$ is a well defined periodic function, that is, 
\begin{equation}
f_p({\bf x} + {\bf u}) \ = \ f_p({\bf x}) \ {\rm for\ all}\ {\bf u} \in {\bf Z}^d
\end{equation}
which can be expressed by a Fourier series of the form
\begin{equation}
f_p({\bf x}) \ = \ \sum_{{\bf u} \in {\bf Z}^d} \ \hat{f}(2 \pi {\bf u}) \ e^{2\pi\ i \ {\bf 1.u}}
\end{equation}
If $f,\ g: {\bf R}^d \rightarrow {\bf R}$ are two absolutely-integrable functions,
the convolution $f * g$ is also absolutely-integrable and is defined as
\begin{equation}
(f*g)({\bf x}) \ = \ \int_{{\bf u} \in {\bf R}^d} \ f({\bf u}) g({\bf x} - {\bf u}) \ d\mu ({\bf u})
\end{equation}
and the fourier transform of $f*g$ is $\hat{f}\hat{g}$.
  
For $r \geq 0$, The set
\begin{equation}
B_{d,\sigma,r}({\bf u})  \ = \ \{{\bf x} \in {\bf R}^d \ : \ \parallel {\bf x} - {\bf u} \parallel_\sigma \leq r \} 
\end{equation}
is termed the {\em $\sigma-$ball} of radius $r$ in ${\bf R}^d$, centered at ${\bf u}$.  
The volume of $B_{d,\sigma,r}({\bf u})$ is denoted by $V_{d,\sigma,r}$.  Clearly, 
\begin{equation}
V_{d,\infty,r} \ = \ (2r)^d
\end{equation}
For $\sigma = 2$  \cite{ref:Weisstein}
\begin{equation}
V_{d,2,r} \ = \ \frac{\pi^{d/2}\ r^d }{\Gamma(1 \ +\  d/2)}
\end{equation}
where $\Gamma$ is the gamma function.  
The surface area of $B_{d,\sigma,r}({\bf u})$ is denoted by 
$A_{d,\sigma,r}$, and it is easy to show that $A_{d,\infty,r} = 2d (2r)^{d-1}$
and that $A_{d,2,r} = d V_{d,2,r}/r$.
In $Q_n^{(\sigma,d)}(r)$, let $\beta_{d,\sigma,r}$ be the probability that
two vertices $i,j$ are connected.  Clearly,
if $0 \leq r \leq 1/2$, $\beta_{d,\sigma,r} =  V_{d,\sigma,r}$.

The Bessel's function of the first kind \cite{ref:Bessel} with parameter $\nu$
is denoted by $J_\nu$.  The following result is well known:
\begin{equation} \label{eq:BallFourierTransform}
\hat\Xi_{B_{d,2,r}(0)}({\bf \omega}) \ = \ (2\pi r )^{d/2}
 \frac{J_{d/2}(r \parallel {\bf\omega} \parallel_2)}{\sqrt{\parallel {\bf \omega} \parallel_2}}
\end{equation}

\section{The probability that a particular labeled cycle appears in $G_n$}

A labeled cycle in $G_n$ of length $q \leq n$ is a sequence of vertices
${\bf y_1}, {\bf y_2}, \ldots {\bf y_q}$ such
that $\{{\bf y_i}, {\bf y_{i+1}}\} \in E(G_n)$ for $i=1,2,\ldots q-1$,
and $\{ {\bf y_q},{\bf y_1}\} \in E(G_n)$.  
Let $\Theta(G_n, q)$ denote the probability
that this labeled cycle is present in $G_n$.  
In both the GR and ER graph, this probability 
does not depend on the particular labeled cycle whose
existence is in question.   Thus,  when $G_n$ is
either an ER or a GR graph, 
\begin{equation} \label{eq:ThetaEq}
\Theta(G_n, q)  \ =  \Theta(G_m,q), \ \ n,m \geq q.
\end{equation}

When $G_n = H(n,p)$, $\Theta(G_n,q)$ can be calculated very easily.
Let $n > 0$ and $1 < q \leq n$.
If $G_n = H(n,p)$, then
the existence of a $q-$cycle in $G_n$ implies the presence
of $q$ edges if $q > 2$, and $q-1$ edges if $q=2$.  In the
ER random graph $H(n,p)$, the presence of an edge is independent
of the presence of the others.  Thus,
\begin{equation} \label{eq:ERCycle}
\Theta(H(n,p),q) = \ \left\{ \begin{array}{ll}  p & {\rm if\ q=2} \\ p^q & {\rm if\ q > 2} \end{array} \right.
\end{equation}

In the case of the geometric random graph $Q_n^{(\sigma,d)}(r)$, things are more complicated
because the edges are not necessarily independent.    Our main result is an
exact characterization of $\Theta(Q_n^{(\sigma,d)}(r)$ for any $\sigma$, $d$.
\begin{theorem} \label{thm:Cycle}
Let $0 < r \leq 1/2$, and $q > 1$.  Then
\begin{equation}
 \Theta(Q_n^{(\sigma,d)}(r),q)  \ =  \ \left\{
\begin{array}{ll} \beta_{d,\sigma,r} & {\rm if }\ q=2 \\
\sum_{{\bf m} \in {\bf Z}^d}  \ 
\hat{\Xi}_{B_{d,\sigma,r}(0)}^q (2 \pi {\bf m}) & {\rm if}\ q>2
\end{array} \right.
\end{equation}
\end{theorem}
\noindent
{\bf Proof:} 
Let ${\bf x}_1, {\bf x}_2, \ldots {\bf x}_q$ be the $q > 1$ random points
which form the labeled cycle of length $q$ (these points are uniformly
distributed in $T_d$).  Then, $\Theta(Q_n^{(\sigma,d)}(r),q)$ is  equal to the probability
that for $i=1,2, \ldots q-1$, 
\begin{equation}
\parallel {\bf x}_i - {\bf x}_{i+1} \parallel_{\sigma} \ \leq \ r
\end{equation}
and $\parallel {\bf x}_q - {\bf x}_1 \parallel_{\sigma} \ \leq \ r$.
Clearly, if $q=2$, then the required probability is just $\beta_{d,\sigma,r}$.

Assume that $q > 2$.  We decompose $\Theta(Q_n^{(\sigma,d)}(r),q)$ as follows:
\begin{eqnarray}
\Theta(Q_n^{(\sigma,d)}(r),q) & = & 
\Pr( \parallel {\bf x}_i - {\bf x}_{i+1} \parallel_\sigma \leq r, 
\ i = 1, 2, \ldots q-1, \ {\rm and} \ \parallel {\bf x}_1 - {\bf x}_q \parallel_\sigma \leq r) \nonumber \\ 
   & = & \Pr (\parallel {\bf x}_1 - {\bf x}_q \parallel_\sigma \leq r  \ /  \ \parallel {\bf x}_i - {\bf x}_{i+1} \parallel_\sigma \leq r, \ i = 1, 2, \ldots q-1)  \nonumber \\
	& & \times \Pr(  \parallel {\bf x}_i - {\bf x}_{i+1} \parallel_\sigma \leq r, \ i = 1, 2, \ldots q-1). \label{eq:3.2.a}
\end{eqnarray}
Clearly, since we are looking at i.i.d. points on the
unit torus $T_1$, the events $\parallel {\bf x}_1 - {\bf x}_{2} \parallel_\sigma \leq r$,  $\parallel {\bf x}_{2} - {\bf x}_{3} \parallel_\sigma \leq r$,
$\ldots$ $\parallel {\bf x}_{q-1} - {\bf x}_q \parallel_\sigma \leq r$ are
independent of each other, and the probability of occurence of each is $\beta_{d,\sigma, r}$.  Hence, 
\begin{equation}
\Pr(  \parallel x_i - x_{i+1} \parallel_\sigma \leq r, \ i = 1, 2, \ldots q-1) \ = \ \beta_{d,\sigma, r} ^{q-1}.
\end{equation}
Thus, we can write
\begin{equation}
 \Theta(Q_n^{(\sigma,d)}(r),q)  \ = \  A_{d,\sigma,q}(r) \times \beta_{d,\sigma, r}^{q-1}  \label{eq:3.2}
\end{equation}
where
\begin{displaymath}
 A_{d,\sigma,q}(r) =   \Pr (\parallel x_1 - x_q \parallel_\sigma \leq r\ / \ \parallel x_i - x_{i+1} \parallel_\sigma \leq r, \ i = 1, 2, \ldots q-1).
\end{displaymath}
We can interpret $A_{d,\sigma,q}(r)$ in the following manner.  
Consider a random walk in ${\bf R}^d$ starting from the origin ${\bf w_1} = {\bf 0}$.  
A point ${\bf u_1}$ is chosen uniformly in the ball 
$B_{d,\sigma,r}({\bf 0})$. The walk then moves to ${\bf w_2} = {\bf w_1} + {\bf u_1}$. 
Continuing in this manner, if the current point is ${\bf w_k}$, 
the walk moves to ${\bf w_{k+1}} = {\bf w_k} + {\bf u_k}$ 
where ${\bf u_k}$ is chosen uniformly in the ball 
$B_{d,\sigma,r}({\bf 0})$.   
Since all points ${\bf m} \in {\bf Z}^d$ map to
the origin ${\bf 0}$ in the unit torus, 
\begin{equation}
A_{d,\sigma,q}(r)  \ = \ \Pr\left( \ {\bf w_q} + {\bf m} \in B_{d,\sigma,r}({\bf 0}) \ {\rm for \ some \ \ m \in {\bf Z}^d }\right)
\end{equation}

Each ${\bf u_i}$ is generated uniformly from $B_{d,\sigma,r}({\bf 0})$, and thus,
the probability density function of each ${\bf u_i}$ is
\begin{equation}
p_u({\bf x}) \ = \ \frac{\Xi_{B_{d,\sigma,r}(0)}({\bf x})}{\beta_{d,\sigma,r}} 
\end{equation}
Then, the probability density function of ${\bf w_k}$
is the $k-1$ fold convolution
\begin{equation}
p_{k}({\bf x}) \ = \ (p_u * p_u * \ldots * p_u)({\bf x}) 
\end{equation}
and the Fourier transform of $p_{k}$ is
\begin{equation} \label{eq:FT}
{\hat{p}_{k}}(\omega) \ = \ \left( 
\frac{\hat{\Xi}_{B_{d,\sigma,r}({\bf 0})}({\bf \omega})}{\beta_{d,\sigma,r}}  \right)^{k-1}
\end{equation}
Define the periodic function $s_r({\bf x})$
as follows
\begin{equation}
s_r({\bf x}) \ = \ \sum_{{\bf m} \in {\bf Z}^d} \ \Xi_{B_{d,\sigma}(r)} ({\bf x} - {\bf m})
\end{equation}
Then, since $r \leq 1/2$,  
\begin{equation} \label{eq:Int1}
A_{d,\sigma,q}(r)  \ = \ \int_{{\bf x} \in {\bf R}^d} \ s_r({\bf x}) p_{q}({\bf x}) d\mu({\bf x})
\end{equation} 
The periodic function $s_r({\bf x})$ has a Fourier series representation 
\begin{equation}
s_r({\bf x})  \ = \   \sum_{{\bf m} \in {\bf Z}^d} \ c_{\bf m} e^{2\pi \  i \  {\bf m . x} }   
\end{equation}
with
\begin{equation}
c_{\bf m} \ = \ \hat{\Xi}_{B_{d,\sigma,r}({\bf 0})}(2 \pi {\bf m})
\end{equation}
Thus, 
\begin{equation} \label{eq:Int2}
A_{d,\sigma,q}(r) \ = \  \int_{{\bf x} \in {\bf R}^d} \ \sum_{{\bf m} \in {\bf Z}^d} c_{\bf m} p_{q}({\bf x}) 
e^{2\pi  i {\bf m.x}} d\mu({\bf x})
\end{equation}
Observe that $p_{q}({\bf x}) = 0$  when $\parallel {\bf x} \parallel_\sigma > qr$.  Thus the integral
in Eq. (\ref{eq:Int2}) can be considered to be 
over a compact set, and the
order of summation and integration can then be exchanged \cite{ref:Young}, and we can write 
\begin{equation} \label{eq:Int3}
A_{d,\sigma,q}(r)  \ = \ \sum_{{\bf m} \in {\bf Z}^d} \ \int_{{\bf x} \in {\bf R}^d}\  c_{\bf m} p_{q}({\bf x}) e^{2\pi i {\bf m. x}} 
d\mu({\bf x})
\end{equation}

For any absolutely integrable $f:{\bf R}^d \rightarrow {\bf R}$, we have
\begin{displaymath}
\int_{{\bf x} \in {\bf R}^d} \ f({\bf x})d\mu ({\bf x}) \ = \ \hat{f}({\bf 0})
\end{displaymath}
Also, by the frequency shift property, 
the Fourier transform of $f(x) e^{i{\bf a.x}}$ is ${\hat{f}}(\omega - {\bf a})$.
Using these facts, we obtain
\begin{equation} \label{eq:AF}
A_{d,\sigma,q}(r)  \ = \   
\sum_{{\bf m} \in {\bf Z}^d}  \ 
\hat{\Xi}_{B_{d,\sigma,r}({\bf 0})}(2 \pi {\bf m})\ \hat{p}_{q}(- 2\pi {\bf m})
\end{equation}
From Eq. (\ref{eq:FT}) and Eq. (\ref{eq:AF}), the theorem follows. $\Box$

Using Theorem \ref{thm:Cycle}, we can obtain series representations for $\Theta$
in terms of the Fourier transform $\hat{\Xi}_{B_{d,\sigma,r}(0)}({\bf \omega})$.
This Fourier transform is relatively easy to compute for $\sigma = \infty$ and for
$\sigma = 2$.

\begin{corollary} \label{cor:RGG}
Let $n > 0$, $0 < r \leq 1/2$, and $1 < q \leq n$.
\begin{equation} \label{eq:RGGCycleInfty}
\Theta( Q_n^{(\infty,d)}(r),q) 
                \ = \ \left\{ \begin{array}{ll} (2r)^d & {\rm if \ q = 2} \\
		\left( 2r \right)^{dq}  \ 
		\left(\ 1  \ + \  2 \sum_{k=1}^{\infty} \left({\rm sinc}(2\pi k r)\right)^q \ \right)^d &
		{\rm if \ q > 2} \end{array} \right. 
\end{equation}
\end{corollary}
{\bf Proof: }
Since $\beta_{d,\infty,r} \ = \ (2r)^d$, the first part of
Eq. (\ref{eq:RGGCycleInfty}) (for $q=2$) follows from Theorem \ref{thm:Cycle}.

Assume that $q > 2$.  Since we are using the $L_{\infty}$ norm, 
each of the $d$ projections of the points
${\bf x}_1, {\bf x}_2, \ldots {\bf x}_q$ must induce a cycle in $T_1$.
Since the projections are independent of each other, it follows
that
\begin{equation} \label{eq:C1E1}
\Theta(Q_n^{(\infty,d)}(r),q) \  = \  \left( \Theta(Q_n^{(\infty,1)}(r),q) \right)^d.
\end{equation}
It is easy to see that
\begin{equation} \label{eq:C1E2}
\hat{\Xi }_{B_{1,\infty,r}(0)}(\omega) \ = \ 2r \ {\rm sinc}(\omega r)
\end{equation}
Using Eq. (\ref{eq:C1E2}) and Eq. (\ref{eq:C1E1}) together with Theorem \ref{thm:Cycle}
we obtain the required expression (we have used  ${\rm sinc}(x) \ = \ {\rm sinc}(-x)$
to rewrite the series). $\Box$

\begin{corollary} \label{cor:RGG2}
Let $n > 0$, $d > 1$, and $1 < q \leq n$. Then
\begin{equation} \label{eq:RGGCycle2}
\Theta( Q_n^{(2,d)}(r),q) 
                \ = \ \left\{ \begin{array}{ll} 
		  V_{d,2,r}  & {\rm if \ q = 2} \\
		  V^{q}_{d,2,r}\ + \  (2\pi r )^{dq/2}\  
		  \sum_{k=1}^{\infty}\ \psi_d(k)\ 
		  \left(\frac{J_{d/2}(2\pi \ r \ \sqrt{k})}{ \left(2\pi\ \sqrt{k}\right)^{1/2}} \right)^q &  {\rm if \ q > 2} \end{array} \right. 
\end{equation}
where $\psi_d(k)$ is the number of solutions ${\bf x} \in {\bf Z}^d$ to the equation
$\parallel {\bf x} \parallel_2 \  = \ k$.
\end{corollary}
\noindent
{\bf Proof:} 
The proof follows immediately from Eq. (\ref{eq:BallFourierTransform})
and Theorem \ref{thm:Cycle}.  $\Box$

\vspace{0.1in}
\noindent
{\bf Remark:}  In order to compute the series in Eq. (\ref{eq:RGGCycle2}),
we need to evaluate the function $\psi_d(k)$.   The following 
recurrence can be used:
\begin{displaymath}
\psi_1(k) \ = \  
\left\{ \begin{array}{ll} 1 & {\rm if}\ k = 0 \\  2 & {\rm if} \ k \neq 0\ {\rm and}\ k = m^2 \ {\rm for\ some} \ m \in {\bf Z} \\
    0 & {\rm otherwise} \end{array} \right.
\end{displaymath}
and if $d > 1$, 
\begin{displaymath}
\psi_d(k) \ = \ \sum_{0 \leq m \leq \sqrt{k}} \psi_{d-1}(k-m^2) 
\end{displaymath}

\section{The expected number of Hamilton cycles in $Q_n^{(2,d)}(r)$}

The Hamilton cycle problem in geometric random graphs
has been studied in \cite{ref:RGGHamilton}, in which
the authors show that the threshold for the existence of
a Hamilton cycle in a geometric random graph (in the
unit cube) is the same as that for 2-connectivity.
The number of Hamilton cycles in a random graph\footnote{The random graph model used in \cite{ref:RGHamiltonCount}
starts with an empty graph on $n$ vertices, and produces a sequence of graphs
by adding new edges with equal probability.  A threshold is then a position in the
sequence at which a property becomes true with high probabililty.} \cite{ref:RGHamiltonCount}
also shows a sharp thresholding property.

Using $\Theta(G_n,n)$, we can directly get the
expected number of Hamilton cycles in $G_n$.
Denote the expected number of
Hamilton cycles in the random graph $G_n$ by $\tau(G_n)$.
For $n > 2$, the number of labeled Hamilton cycles in a complete graph on $n$ vertices
is $(n-1)!/2$.   
It follows that, for $n > 2$,
\begin{equation}
\tau(G_n) \ = \  \Theta(G_n,n) \ (n-1)!/2  
\end{equation}
because the probability of each such labeled cycle being present
is $\Theta(G_n,n)$.

\begin{figure}
\begin{centering}
\centerline{\includegraphics[width=3.0in,height=2.75in]{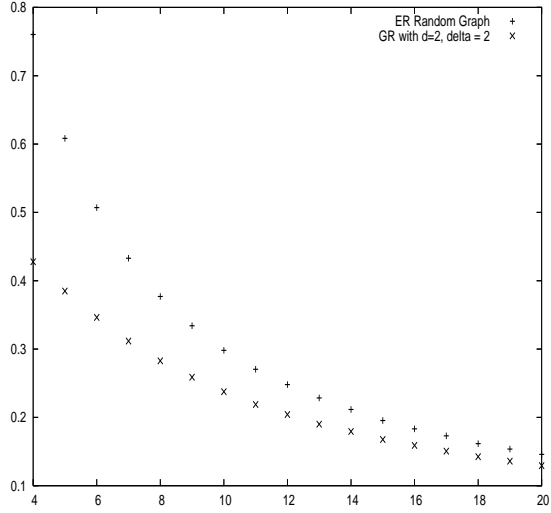}}
\caption{Threshold for $\tau(G_n) \geq 1$  plotted
as a function of the $n$ for the ER graph and for
the GR graph with $d=2$, $\sigma = 2$}
 \label{fig:HCThres}
\end{centering}
\end{figure}

Consider the threshold for $G_n$ defined as the
smallest edge-probability such that $\tau(G_n) \geq 1$.
We can use Corollaries \ref{cor:RGG} and \ref{cor:RGG2} to compute
this threshold when $G_n = Q_{n}^{(2,d)}(r)$ and $G_n = Q_n^{\infty,d}(r)$,
and contrast this threshold with that for the ER graph $H(n,p)$.
In Figure \ref{fig:HCThres}, we show the thresholds
obtained for $H(n,p)$ and $Q_n^{(2,\sigma)}(r)$.
The computed threshold for
the geometric random graph is lower than that for
the ER graph.  However, the difference between
the two thresholds reduces as $n$ increases. 
Asymptotically, the 
threshold for the appearance of a Hamilton cycle
seems to be similar in the GR graph and the 
ER random graph (this threshold is of the order
$\log(n) /n$ \cite{ref:RGGHamilton}).  An explanation
for this is that as $n$ increases, the end points of
a path of length $n$ become less correlated (recall the
random walk argument used in the proof of Theorem \ref{thm:Cycle}), 
and thus, the probability of an edge between the end points
of the path is close to the edge probability. 

\section{The expected value of the determinant and the permanent of $A_{G_n}$}

Let $F_{G_n}(x)$ be the
matrix $xI + A_{G_n}$.  Define the two polynomials
\begin{equation}
\Lambda_{G_n}(x) \ = \ \det (F_{G_n}(x))), 
\end{equation}
and 
\begin{equation}
\Gamma_{G_n}(x) \ = \ {\rm per} (F_{G_n}(x)).
\end{equation}
The polynomials $\Lambda_{G_n}(x)$ and $\Gamma_{G_n}(x)$ have 
coefficients which are random variables.
In particular, the coefficients in $\Lambda_{G_n}$
are symmetric functions of the 
eigenvalues of $A_{G_n}$.
Define 
\begin{equation}
\bar{\Lambda}_{G_n} (x) \ = \ E(\Lambda_{G_n}(x))
\end{equation}
and 
\begin{equation}
\bar{\Gamma}_{G_n} (x) \ = \ E(\Gamma_{G_n}(x))
\end{equation}
where the expectation of a polynomial $p(x)$ is the polynomial $\bar{p}(x)$ whose 
coefficients are the expectations of the corresponding coefficients in $p(x)$.

The coefficient of $x^k$ in $\bar{\Lambda}_{G_n}(x)$ is the 
expected value of the elementary symmetric function of degree $n-k$ evaluated
at the eigenvalues of $A_{G_n}$.  In particular, the constant term in
$\bar{\Lambda}_{G_n}(x)$ is the expected value of the determinant of
$A_{G_n}$, so that the expected value of
the determinant of $A_{G_n}$ is $\bar{\Lambda}_{G_n}(0)$.   
The coefficient of $x^k$ in $\bar{\Gamma}_{G_n}(x)$ is
the expected number of cycle covers across all subgraphs of $G_n$
with $n-k$ vertices.  Also, the expected value
of the permanent of $G_n$ is $\bar{\Gamma}_{G_n}(0)$. 

There is a strong connection between cycles and permutations, and
between permutations and determinants (and permanents). We 
expect that the characterization of $\Theta(G_n,q)$ will help
determine the behaviour of the determinant (and permanent).
More concretely, we show that
\begin{theorem} \label{thm:MainThm}
Let $G_n$ be a random graph on $n > 0$ vertices ($G_n$ is either the ER graph
or the GR graph).   Then, 
for $n \geq 1$, the polynomials $\bar{\Lambda}_{G_n}(x)$ and
$\bar{\Gamma}_{G_n}(x)$ satisfy the recurrence relations
\begin{equation} 
\bar{\Lambda}_{G_n} (x) \ =  \ x \bar{\Lambda}_{G_{n-1}}(x) \ +  \ 
\sum_{q=2}^n \ (-1)^{q-1} \frac{n-1!}{n-q!} \ 
\  \Theta(G_n,q) \ \bar{\Lambda}_{G_{n-q}}(x) 
\end{equation}
and 
\begin{equation}
\bar{\Gamma}_{G_n} (x) \ =  \ x \bar{\Gamma}_{G_{n-1}}(x) \ +  \ 
\sum_{q=2}^n \ \frac{n-1!}{n-q!}  
\ \Theta(G_n,q) \ \bar{\Gamma}_{G_{n-q}}(x) 
\end{equation}
with initial conditions $\bar{\Lambda}_{G_0}(x) = \bar{\Gamma}_{G_0}(x) = 1$.
\end{theorem}

\noindent
{\bf Proof:}
We start with the following formulas for the determinant
and the permanent.  If $B = [b_{ij}]$ is an $n \times n$
matrix, then 
\begin{equation}
\det(B) \ = \  \sum_{\sigma \in S_n} (-1)^{sign({\sigma})} \prod_{i=1}^{n} b_{i\sigma (i)}
\end{equation}
and 
\begin{equation}
{\rm per}(B) \ = \  \sum_{\sigma \in S_n}  \prod_{i=1}^{n} b_{i\sigma (i)}
\end{equation}
where $S_n$ is the group of permutations of $\{1,2, \ldots n\}$.

Each permutation $\sigma \in S_n$ can be uniquely
decomposed into a set  of 
disjoint cycles on $\{ 1, 2, \ldots n \}$.  
Each cycle $C$ in the disjoint cycle-decomposition of a permutation
is of the form $(i_1 i_2 \ldots i_q)$, where
$\sigma(i_r) \ = \ i_{r+1}, r = 1,2, \ldots q-1$
and $\sigma(i_q) \ = i_1$.  The sign of the
cycle $C$ is $sign(C) = (-1)^{|C|-1}$, where $|C|$
is the number of elements in $C$. 
The sign of the permutation is then the product of signs
of the cycles into which $\sigma$ is decomposed.  We will say that the
pair $(i,j) \in C$ if $i, j$ are consecutive elements
in the cycle $C$ ($i_0$ is considered to be after $i_q$).
Then, given $\sigma$, we have
\begin{equation}
 \prod_{i=1}^{n} b_{i\sigma (i)} \ = \ \prod_{C \in \sigma} \prod_{(i,j) \in C}\  b_{ij}
\end{equation}
For a cycle $C$, We define 
\begin{equation}
w_B(C) \ = \ \prod_{(i,j) \in C} b_{ij}
\end{equation}
Then, 
\begin{equation}
\det(B) \ = \  \sum_{\sigma \in S_n} \prod_{C \in \sigma} \ (-1)^{|C|-1} w_B(C)
\end{equation}
and 
\begin{equation}
{\rm per} (B) \ = \  \sum_{\sigma \in S_n} \prod_{C \in \sigma}  \ w_B(C).
\end{equation}

Let $B = F_{G_n}(x)$. For a cycle  $C = (i_1 i_2 \ldots i_q)$ in some permutation,
we see that if $q > 1$, then 
\begin{equation}
E(w_B(C)) \ = \ \Theta(G_n,q)
\end{equation} 
and if $q = 1$, then 
\begin{equation}
E(w_B(C)) \ = \ x.
\end{equation}
For convenience, we set $\Theta(G_n,1) = x$.

Also, if $C_1, C_2, \ldots C_t$ are vertex-disjoint cycles in 
$G_n$, then the presence of $C_i$ is independent of the presence
of $C_j$ for $j \neq i$, and 
\begin{equation}
E(\prod_{i=1}^t \ w_B(C_i)) \ = \ \prod_{i=1}^t \ E(w_B(C_i)).
\end{equation}
It follows that 
\begin{equation}
\bar{\Lambda}_{G_n}(x) \ = \  \sum_{\sigma \in S_n} \prod_{C \in \sigma} (-1)^{|C|-1} \Theta({G_n},|C|).
\end{equation}
Similarly, 
\begin{equation}
\bar{\Gamma}_{G_n}(x) \ = \  \sum_{\sigma \in S_n} \prod_{C \in \sigma} \Theta({G_n},|C|)
\end{equation}
The counting of permutations $\sigma \in S_n$ can be carried out
by fixing the cycle $C$ which contains $1$ and counting permutations
of elements not in $C$.  For $1 \leq q \leq n$, Let $D_q$ be the set of cycles of length $q$ which
contain $1$.
We observe that 
\begin{displaymath}
| D_q | \ = \ (q-1)! \ \left( \begin{array}{c} n-1 \\ q-1 \end{array} \right),
\end{displaymath}
because each cycle in $D_q$ is determined by the choice of $q-1$ elements
(other than $1$) out of $n-1$ elements, and there are $q-1!$ 
distinct cycles on $q$ elements.

Let ${\bf N} = \{1,2, \ldots, n \}$ and let $P(A)$ be
the set of permutations of the set $A \subset {\bf N}$. 
Then, we can write
\begin{equation} \label{eq:Rec0}
\sum_{\sigma \in S_n} \prod_{C \in \sigma} (-1)^{|C|-1}\ \Theta({G_n},|C|)
\end{equation}
as 
\begin{equation} \label{eq:Rec1}
\sum_{q=1}^{n} \left( \sum_{C \in D_q} (-1)^{|C|-1}\ \Theta({G_n},q) \nonumber 
\left( \sum_{\sigma \in P({\bf N} - C)} \prod_{D \in \sigma} (-1)^{|D|-1} \ \Theta({G_n},|D|) \right)\right)
\end{equation}
where the innermost summation over $P(A)$ is taken to be $1$ if $A = \phi$.
Since $|C| = q$ for each $C \in D_q$, we can rewrite Eq. (\ref{eq:Rec1}) (using
Eq. (\ref{eq:ThetaEq}) to replace $\Theta({G_n},|D|)$ by $\Theta (G_{n-q},|D|)$) 
as
\begin{equation} \label{eq:Rec2}
\sum_{q=1}^{n} \left( \begin{array}{c} n-1 \\ q-1 \end{array} \right) \ (q-1)!\ (-1)^{q-1}\ \Theta({G_n},q) \nonumber 
\left( \sum_{\sigma \in P({\bf N} - C)} \prod_{D \in \sigma} (-1)^{|D|-1} \ \Theta({G_{n-q}},|D|) \right).
\end{equation}
The inner summation in Eq. (\ref{eq:Rec2}) is just $\bar{\Lambda}_{G_{n-q}}(x)$, and thus,
the recurrence relation for $\bar{\Lambda}_{G_n}(x)$ follows. 
The recurrence relation for $\bar{\Gamma}_{G_n}(x)$ can be shown to hold in 
a similar manner, completing the proof of Theorem \ref{thm:MainThm}.  $\Box$

\vspace{0.1in}
\noindent
{\bf Remark:}  The result in Theorem \ref{thm:MainThm} holds for any
random graph $G_n$ in which 
the probability of appearance of a labeled
cycle depends only on its length and  
the probability of appearance of a set of vertex-disjoint cycles is
the product of probabilities of appearance of the elements in this set.

\vspace{0.1in}

For $n > 0$, $0 < k \leq n$, let $F_{n,k}(t_1, t_2, \ldots t_n)$ denote the elementary symmetric function
\begin{equation}
F_{n,k} (t_1, t_2, \ldots t_n) \ = \ \sum_{\{i_1,i_2, \ldots i_k \} \in \{1,2,\ldots,n\}} \ t_{i_1} t_{i_2} \ldots t_{i_k}
\end{equation}
For $k=0$, define $F_{n,k} = 1$, and define $F_{n,k} = 0$ if $n < k$ or if $k < 0$.
Now, let $\hat{F}_{n,k}$ denote the expected value of $F_{n,k}$ evaluated
on the $n$ eigenvalues of $A_{G_n}$.  Then, the expected value of the determinant of $A_{G_n}$
is just $\hat{F}_{n,n}$.  Then, we have the following corollary of Theorem \ref{thm:MainThm}.
\begin{corollary}
For the random graph $G_n$, if $n > 0$, and $0 < k \leq n$, then
\begin{equation}
\hat{F}_{n,k} \ = \ \hat{F}_{n-1,k} \ + \ \sum_{q=2}^{n} \ (-1)^q \frac{n-1!}{n-q!} \ \Theta(G_n,q) \ \hat{F}_{n-q,k-q}
\end{equation}
\end{corollary}
\noindent
{\bf Proof:} Follows from Theorem \ref{thm:MainThm} by noting that the coefficient of 
$x^k$ in $\bar{\Lambda}$ is $\hat{F}_{n,n-k}$. $\Box$

\vspace{0.2in}
Note that in both models, if the edge probability is $1$, 
then $\Theta(G_n,q) = 1$, and $G_n$ is always the
complete graph, so that the expected value of the determinant
of $G_n$ is $(-1)^{n-1} \times (n-1)$.
Using Theorem \ref{thm:MainThm}, 
we obtain the following identity for $n > 0$:
\begin{equation}
n \ = \ 1 \ + \  \sum_{q=2}^{n}  \ \frac{n-1 !}{n-q!} \times ((n-q)-1)  
\end{equation} 
Also, the permanent of the complete
graph on $n$ vertices is the number of
derangements of the set ${\bf N} = \{ 1,2, \ldots n \}$.   Thus, 
the recurrence proved in Theorem \ref{thm:MainThm}
yields the following identity for the number of
derangements $d_n$ of {\bf N}
\begin{equation}
d_n \ = \ \sum_{q=2}^{n} \frac{n-1 !}{n-q !} \ d_{n-q}, 
\end{equation} 
with the initial conditions $d_1 = 0$, and $d_0 = 1$.

We use these ecurrence relations to compute 
these expected values for $n \leq 20$ in the GR and ER
models\footnote{The recurrence relations were directly
computed using {\em long double} precision arithmetic. For  
higher values of $n$ one would need to use higher precision arithmetic.}.  
Some interesting conclusions can be 
drawn from these calculations.

\begin{figure}
\begin{centering}
\centerline{\includegraphics[width=3.0in,height=2.75in]{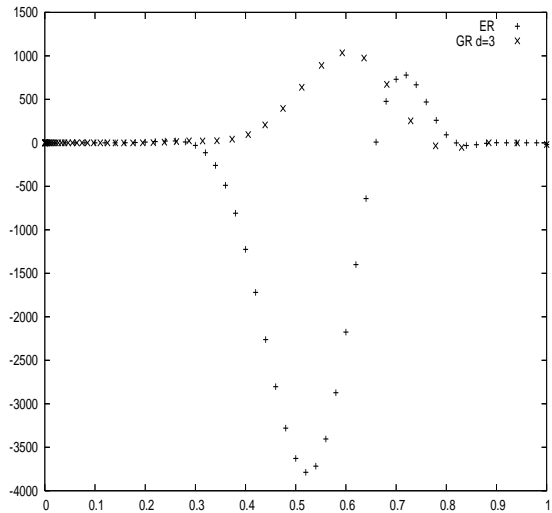}}
\caption{The expected value of the determinant plotted as a function of the
edge probability for $n=20$ in the ER and GR (with $d=3$, $\sigma=\infty$) models.}
 \label{fig:DetPlotx20xERxGR3}
\end{centering}
\end{figure}

Consider the plot in Figure \ref{fig:DetPlotx20xERxGR3}, in which 
we compare the behaviour of the determinant of
$G_{20}$ as a function of the edge probability.
The graph has been plotted for the ER graph and 
for the GR graph with $d=3$.  The
behaviour of the determinant in the two models
is quite different, and clearly, so is the
distribution of $G_n$.  

\begin{figure}
\begin{centering}
\centerline{\includegraphics[width=3.0in,height=2.75in]{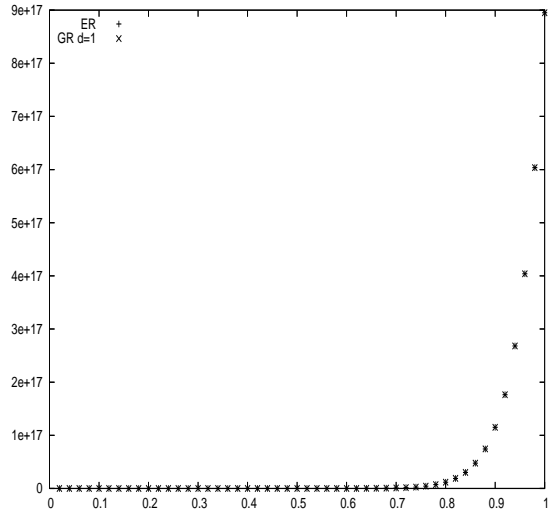}}
\caption{The expected value of the permanent plotted as a function of the
edge probability for $n=20$ in the ER and GR (with $d=1$, $\sigma=\infty$) models.}
 \label{fig:PerPlotx20x1}
\end{centering}
\end{figure}

\begin{figure}
\begin{centering}
\centerline{\includegraphics[width=3.0in,height=2.75in]{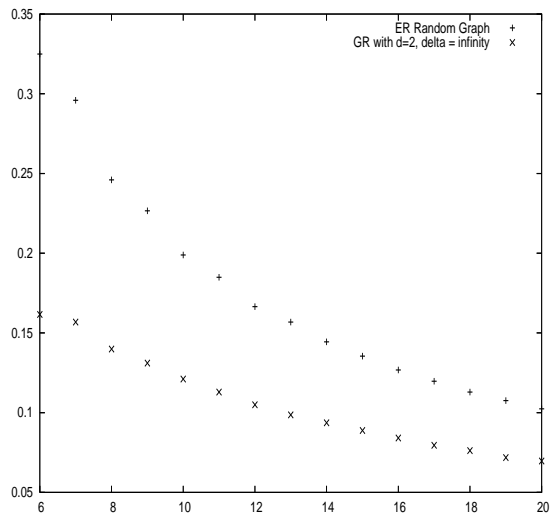}}
\caption{Threshold for the expected value of the permanent  plotted
as a function of the $n$ for the ER graph and for
the GR graph with $d=2$, $\sigma = \infty$}
 \label{fig:PERThres}
\end{centering}
\end{figure}

In Figure \ref{fig:PerPlotx20x1}, we show a plot of
the expected value of the permanent of $A_{G_n}$ (for $n=20$)
as a function of the edge probability in the ER and GR ($d=1$) 
models.  
We can also define a threshold
for the expected permanent as the smallest edge probability
for which the expected value of the permanent is $\geq 1$.
A comparison of this threshold for the GR and ER graphs shows
that this threshold is lower for the GR graph, but the two
thresholds come closer as $n$ increases (see Figure \ref{fig:PERThres}).
Thus, the permanent of the GR graph grows more rapidly than
that of the ER graph.  This is expected since a labeled cycle
is more likely in the GR graph.


\subsection{Graphs with large determinants}

Looking at Figure \ref{fig:DetPlotx20xERxGR3}, 
we see that for intermediate values of the edge probability, large magnitudes
appear in the plots of the expected value of the
determinant.  
For instance, we observe that, 
in the ER random graph with $n=20$, the largest absolute
value of the determinant is $3787.81$, and this provides a lower
bound on the maximal determinant of a symmetric
$20\times 20$ $0/1$ matrix .

For a general (non-symmetric) $n\times n$
$0/1$ matrix, the determinant is
bounded above by $(n+1)^{(n+1)/2}/2^n$ \cite{ref:Faddeev}.
The number of (possibly non-symmetric)
$n \times n$  $0/1$ matrices which  achieve
this bound  is also known for $n \leq 9$ \cite{ref:Miodrag}.  
However, similar characterizations of
the determinants of {\em symmetric} $0/1$ matrices
are not so common.  For example,  in \cite{ref:Fallat}, 
the authors show that for $n \geq 7$, the maximal determinant of
the adjacency matrix of a $(n-3)$-regular graph on $n$ vertices
is $(n-3) 3^{[n/4] - 1}$.  For $n=20$, this works out
to be $1377$ which is less than the largest observed
determinant value in the evolution of $H(20,p)$.

Thus, the recurrence formula for the
expected value of the determinant seems 
to provide some useful information about the maximal
determinant of a class of 
symmetric $0/1$ matrices (in effect, we have 
a lower bound on the largest value of such
determinants).  Also, if
the expected determinant is large, then
it may be possible to find
a symmetric $0/1$ matrix with 
large determinant by using a Monte Carlo sampling
approach.  An estimate
of the second moment of the determinant of
the random graph will throw more light on
this possibility.

\section{Conclusions}

We have derived an exact characterization of the probability
of existence of a labeled cycle in geometric random 
graphs on a unit torus with an arbitrary number of dimensions,
and with an arbitrary $L_\sigma$ distance metric).  
This cycle probability can be calculated
in terms of the Fourier transform of the indicator function 
of a ball in $L_\sigma$.  Explicit expressions for this
Fourier transform can be easily computed in the $\sigma = \infty$
and $\sigma = 2$ case.  

From the cycle probability, one gets the expected number of
Hamilton cycles in the geometric random graph.  These exact
expressions complement the asymptotic threshold results
for the existence of Hamilton cycles in geometric random 
graphs (as in \cite{ref:RGGHamilton}).  We observe that
as the edge probability increases, a Hamilton cycle
appears earlier in the GR graph than
in the ER graph.  

The cycle probabilities can also be used to find the
expected values of the determinant (and more generally, the
expected values of the elementary symmetric functions evaluated
at the eigenvalues of the adjacency matrix) and the permanent
of the adjacency matrix of the random graph.  We obtain recurrence relations
for these quantities and illustrate them by a few calculations.
In particular, the determinant exhibits
very different behaviour in the two models.  Also, large
magnitudes of the determinant are observed in the evolution
of the random graphs.  This throws some light on the
as yet unresolved question of the maximal determinant of symmetric 0/1 matrices.


\end{document}